%Paper: q-alg/9512010
%From: chered@math.unc.edu (Ivan Cherednik)
%Date: Thu, 7 Dec 1995 13:07:25 +0500

% formatam.tex -- AMSTeX template file
% Version: June 10, 1992

\input amstex

%\documentstyle{amams} % input Annals of Mathematics macros.
%%%%%%%%%%%%%%%%%%%%%%%%%%%%%%%%%%%%%%%%%%%%%%%%%%%%%%%%%%%
%% amams.sty: AMSTeX Macros for Articles to be published in
%%
%% Annals of Mathematics
%%
%% Princeton University and the
%% Institute for Advanced Study
%%
%% Published by Princeton University Press
%%%%%%%%%%%%%%%%%%%%%%%%%%%%%%%%%%%%%%%%%%%%%%%%%%%%%%%%%%%

%%%%%%%%%%%%%%%%%%%%%%%%%%%%%%%%%%%%%%%%%%%%%%%%%%%%%%%%%%%
%% Variations on AMSPPT.sty written by Amy Hendrickson
%% TeXnology Inc, Brookline, MA
%% 617 738-8029, amyh@ai.mit.edu
%%%%%%%%%%%%%%%%%%%%%%%%%%%%%%%%%%%%%%%%%%%%%%%%%%%%%%%%%%%

\def\spaces{\space\space\space\space\space\space\space\space\space\space}
\def\spacess{\message{\spaces\spaces\spaces\spaces\spaces\spaces\spaces}}
\spacess
\spacess
\message{Annals of Mathematics Style: Current Version: 1.1. June 10, 1992}
\spacess
\spacess
%%%%%%%%%%%%%%%%%%%%%%%%%%%%%%%%%%%%%%%%%%%%%%%%%%%%%%%%

\catcode`\@=11

\hyphenation{acad-e-my acad-e-mies af-ter-thought anom-aly anom-alies
an-ti-deriv-a-tive an-tin-o-my an-tin-o-mies apoth-e-o-ses apoth-e-o-sis
ap-pen-dix ar-che-typ-al as-sign-a-ble as-sist-ant-ship as-ymp-tot-ic
asyn-chro-nous at-trib-uted at-trib-ut-able bank-rupt bank-rupt-cy
bi-dif-fer-en-tial blue-print busier busiest cat-a-stroph-ic
cat-a-stroph-i-cally con-gress cross-hatched data-base de-fin-i-tive
de-riv-a-tive dis-trib-ute dri-ver dri-vers eco-nom-ics econ-o-mist
elit-ist equi-vari-ant ex-quis-ite ex-tra-or-di-nary flow-chart
for-mi-da-ble forth-right friv-o-lous ge-o-des-ic ge-o-det-ic geo-met-ric
griev-ance griev-ous griev-ous-ly hexa-dec-i-mal ho-lo-no-my ho-mo-thetic
ideals idio-syn-crasy in-fin-ite-ly in-fin-i-tes-i-mal ir-rev-o-ca-ble
key-stroke lam-en-ta-ble light-weight mal-a-prop-ism man-u-script
mar-gin-al meta-bol-ic me-tab-o-lism meta-lan-guage me-trop-o-lis
met-ro-pol-i-tan mi-nut-est mol-e-cule mono-chrome mono-pole mo-nop-oly
mono-spline mo-not-o-nous mul-ti-fac-eted mul-ti-plic-able non-euclid-ean
non-iso-mor-phic non-smooth par-a-digm par-a-bol-ic pa-rab-o-loid
pa-ram-e-trize para-mount pen-ta-gon phe-nom-e-non post-script pre-am-ble
pro-ce-dur-al pro-hib-i-tive pro-hib-i-tive-ly pseu-do-dif-fer-en-tial
pseu-do-fi-nite pseu-do-nym qua-drat-ics quad-ra-ture qua-si-smooth
qua-si-sta-tion-ary qua-si-tri-an-gu-lar quin-tes-sence quin-tes-sen-tial
re-arrange-ment rec-tan-gle ret-ri-bu-tion retro-fit retro-fit-ted
right-eous right-eous-ness ro-bot ro-bot-ics sched-ul-ing se-mes-ter
semi-def-i-nite semi-ho-mo-thet-ic set-up se-vere-ly side-step sov-er-eign
spe-cious spher-oid spher-oid-al star-tling star-tling-ly
sta-tis-tics sto-chas-tic straight-est strange-ness strat-a-gem strong-hold
sum-ma-ble symp-to-matic syn-chro-nous topo-graph-i-cal tra-vers-a-ble
tra-ver-sal tra-ver-sals treach-ery turn-around un-at-tached un-err-ing-ly
white-space wide-spread wing-spread wretch-ed wretch-ed-ly Brown-ian
Eng-lish Euler-ian Feb-ru-ary Gauss-ian Grothen-dieck Hamil-ton-ian
Her-mit-ian Jan-u-ary Japan-ese Kor-te-weg Le-gendre Lip-schitz
Lip-schitz-ian Mar-kov-ian Noe-ther-ian No-vem-ber Rie-mann-ian
Schwarz-schild Sep-tem-ber Za-mo-lod-chi-kov Knizh-nik quan-tum Op-dam
Mac-do-nald Ca-lo-ge-ro Su-ther-land Mo-ser Ol-sha-net-sky  Pe-re-lo-mov
in-de-pen-dent ope-ra-tors
}

\Invalid@\nofrills
\Invalid@\usualspace
\newif\ifnofrills@
\def\nofrills@#1#2{\relaxnext@
  \DN@{\ifx\next\nofrills
    \nofrills@true\let#2\relax\DN@\nofrills{\nextii@}%
  \else
    \nofrills@false\def#2{#1}\let\next@\nextii@\fi
\next@}}
\def\usualspace@#1{\ifnofrills@\def\usualspace{#1}\fi}
\def\addto#1#2{\csname \expandafter\eat@\string#1@\endcsname
  \expandafter{\the\csname \expandafter\eat@\string#1@\endcsname#2}}
\newdimen\bigsize@
\def\big@#1#2{{\hbox{$\left#2\vcenter to#1\bigsize@{}%
  \right.\nulldelimiterspace\z@\m@th$}}}
\def\big{\big@\@ne}
\def\Big{\big@{1.5}}
\def\bigg{\big@\tw@}
\def\Bigg{\big@{2.5}}
\def\raggedcenter@{\leftskip\z@ plus.4\hsize \rightskip\leftskip
 \parfillskip\z@ \parindent\z@ \spaceskip.3333em \xspaceskip.5em
 \pretolerance9999\tolerance9999 \exhyphenpenalty\@M
 \hyphenpenalty\@M \let\\\linebreak}
\def\upperspecialchars{\def\ss{SS}\let\i=I\let\j=J\let\ae\AE\let\oe\OE
  \let\o\O\let\aa\AA\let\l\L}
\def\uppercasetext@#1{%
  {\spaceskip1.2\fontdimen2\the\font plus1.2\fontdimen3\the\font
   \upperspecialchars\uctext@#1$\m@th\aftergroup\eat@$}}
\def\uctext@#1$#2${\endash@#1-\endash@$#2$\uctext@}
\def\endash@#1-#2\endash@{%
\uppercase{#1}\if\notempty{#2}--\endash@#2\endash@\fi}
\def\runaway@#1{\DN@{#1}\ifx\envir@\next@
  \Err@{You seem to have a missing or misspelled \string\end#1 ...}%
  \let\envir@\empty\fi}
\newif\iftemp@
\def\notempty#1{TT\fi\def\test@{#1}\ifx\test@\empty\temp@false
  \else\temp@true\fi \iftemp@}

%\comment%%% remove
\font@\tensmc=cmcsc10
\font@\sevenex=cmex7
\font@\sevenit=cmti7
\font@\eightrm=cmr8 % preloaded in plain.tex
\font@\sixrm=cmr6 % preloaded in plain.tex
\font@\eighti=cmmi8     \skewchar\eighti='177 % preloaded
\font@\sixi=cmmi6       \skewchar\sixi='177   % preloaded
\font@\eightsy=cmsy8    \skewchar\eightsy='60 % preloaded
\font@\sixsy=cmsy6      \skewchar\sixsy='60   % preloaded
\font@\eightex=cmex8 %
\font@\eightbf=cmbx8 % preloaded in plain.tex
\font@\sixbf=cmbx6   % preloaded in plain.tex
\font@\eightit=cmti8 % preloaded in plain.tex
\font@\eightsl=cmsl8 % preloaded in plain.tex
\font@\eightsmc=cmcsc10
\font@\eighttt=cmtt8 % preloaded in plain.tex
%\font@\ninerm=cmr9
%\font@\ninei=cmmi9    \skewchar\ninei='177
%\font@\ninesy=cmsy9   \skewchar\ninesy='60
%\font@\nineex=cmex9
%\font@\ninebf=cmbx9
%\font@\nineit=cmti9
%\font@\ninesl=cmsl9
%\font@\ninesmc=cmcsc9
%\font@\ninemsa=msam9
%\font@\ninemsb=msbm9
%\font@\nineeufm=eufm9
%\endcomment%%%

\loadmsam
\loadmsbm
\loadeufm
\UseAMSsymbols

\def\penaltyandskip@#1#2{\relax\ifdim\lastskip<#2\relax\removelastskip
      \ifnum#1=\z@\else\penalty@#1\relax\fi\vskip#2%
  \else\ifnum#1=\z@\else\penalty@#1\relax\fi\fi}
\def\nobreak{\penalty\@M
  \ifvmode\def\penalty@{\let\penalty@\penalty\count@@@}%
  \everypar{\let\penalty@\penalty\everypar{}}\fi}
\let\penalty@\penalty

\def\block{\RIfMIfI@\nondmatherr@\block\fi
       \else\ifvmode\vskip\abovedisplayskip\noindent\fi
        $$\def\endblock{\par\egroup$$}\fi
  \vbox\bgroup\advance\hsize-2\indenti\noindent}
\def\endblock{\par\egroup}

\def\logo@{\baselineskip2pc \hbox to\hsize{\hfil\eightpoint Typeset by
 \AmSTeX}}

%%%%%%%%%%%%%%%%%%%%%%%%%%%%%%%%%%%%%%%%%%%%%%%%%%%%%%%%%%%%%%%
%% Macros for Annals of Mathematics written by Amy Hendrickson
%% TeXnology Inc, Brookline, MA
%% 617 738-8029, amyh@ai.mit.edu
%%%%%%%%%%%%%%%%%%%%%%%%%%%%%%%%%%%%%%%%%%%%%%%%%%%%%%%%%%%%%%%

%% This file includes:
%% 1) Font declarations,
%% 2) Page set up,
%% 3) Title page
%% 4) Section heads,
%% 5) Equation macros, autonumbering equations, etc.,
%% 6) Figure and Table Captions,
%% 7) End matter macros: Bibliography, Appendix, etc.,
%% 8) Footnotes,
%% 9) Theorem type environments
%% 10) Cross-referencing
%% 11) Listing
%% 12) Article and Journal Table of Contents

%%%%%%%%%%%%%%%%%%%%%%%%%%%%%%%%%%%
%% 1) Font declarations,
% Computer Modern fonts

% Small Caps
\font\elevensc=cmcsc10 scaled\magstephalf
\font\tensc=cmcsc10

\font\eightsc=cmcsc10 scaled800

\font\elevenrm=cmr10 scaled \magstephalf%!!!
\font\ninerm=cmr9
\font\eightrm=cmr8
\font\sixrm=cmr6
\font\fiverm=cmr5

\font\eleveni=cmmi10 scaled\magstephalf
\font\ninei=cmmi9
\font\eighti=cmmi8
\font\sixi=cmmi6
\font\fivei=cmmi5
\skewchar\ninei='177 \skewchar\eighti='177 \skewchar\sixi='177
\skewchar\eleveni='177

\font\elevensy=cmsy10 scaled\magstephalf
\font\ninesy=cmsy9
\font\eightsy=cmsy8
\font\sixsy=cmsy6
\font\fivesy=cmsy5
\skewchar\ninesy='60 \skewchar\eightsy='60 \skewchar\sixsy='60
\skewchar\elevensy'60

\font\eighteenbf=cmbx10 scaled\magstep3

\font\twelvebf=cmbx10 scaled \magstep1
\font\elevenbf=cmbx10 scaled \magstephalf
\font\tenbf=cmbx10
\font\ninebf=cmbx9
\font\eightbf=cmbx8
\font\sixbf=cmbx6
\font\fivebf=cmbx5

\font\elevenit=cmti10 scaled\magstephalf
\font\nineit=cmti9
\font\eightit=cmti8

% Fonts for bold math
\font\eighteenmib=cmmib10 scaled \magstep3
\font\twelvemib=cmmib10 scaled \magstep1
\font\elevenmib=cmmib10 scaled\magstephalf
\font\tenmib=cmmib10
\font\eightmib=cmmib10 scaled 800
\font\sixmib=cmmib10 scaled 600

\font\eighteensyb=cmbsy10 scaled \magstep3
\font\twelvesyb=cmbsy10 scaled \magstep1
\font\elevensyb=cmbsy10 scaled \magstephalf
\font\tensyb=cmbsy10
\font\eightsyb=cmbsy10 scaled 800
\font\sixsyb=cmbsy10 scaled 600

\font\elevenex=cmex10 scaled \magstephalf
\font\tenex=cmex10
\font\eighteenex=cmex10 scaled \magstep3

%%%%%%%%%%%%%%%%%%%%%%%%%%%%
%% Font families

\def\elevenpoint{\def\rm{\fam0\elevenrm}%
  \textfont0=\elevenrm \scriptfont0=\eightrm \scriptscriptfont0=\sixrm
  \textfont1=\eleveni \scriptfont1=\eighti \scriptscriptfont1=\sixi
  \textfont2=\elevensy \scriptfont2=\eightsy \scriptscriptfont2=\sixsy
  \textfont3=\elevenex \scriptfont3=\tenex \scriptscriptfont3=\tenex
  \def\bf{\fam\bffam\elevenbf}%
  \def\it{\fam\itfam\elevenit}%
  \textfont\bffam=\elevenbf \scriptfont\bffam=\eightbf
   \scriptscriptfont\bffam=\sixbf
\normalbaselineskip=13.95pt
  \setbox\strutbox=\hbox{\vrule height9.5pt depth4.4pt width0pt\relax}%
  \normalbaselines\rm}

\elevenpoint %%% default fonts and baselineskip

\def\ninepoint{\def\rm{\fam0\ninerm}%
  \textfont0=\ninerm \scriptfont0=\sixrm \scriptscriptfont0=\fiverm
  \textfont1=\ninei \scriptfont1=\sixi \scriptscriptfont1=\fivei
  \textfont2=\ninesy \scriptfont2=\sixsy \scriptscriptfont2=\fivesy
  \textfont3=\tenex \scriptfont3=\tenex \scriptscriptfont3=\tenex
  \def\it{\fam\itfam\nineit}%
  \textfont\itfam=\nineit
  \def\bf{\fam\bffam\ninebf}%
  \textfont\bffam=\ninebf \scriptfont\bffam=\sixbf
   \scriptscriptfont\bffam=\fivebf
\normalbaselineskip=11pt
  \setbox\strutbox=\hbox{\vrule height8pt depth3pt width0pt\relax}%
  \normalbaselines\rm}

\def\eightpoint{\def\rm{\fam0\eightrm}%
  \textfont0=\eightrm \scriptfont0=\sixrm \scriptscriptfont0=\fiverm
  \textfont1=\eighti \scriptfont1=\sixi \scriptscriptfont1=\fivei
  \textfont2=\eightsy \scriptfont2=\sixsy \scriptscriptfont2=\fivesy
  \textfont3=\tenex \scriptfont3=\tenex \scriptscriptfont3=\tenex
  \def\it{\fam\itfam\eightit}%
  \textfont\itfam=\eightit
  \def\bf{\fam\bffam\eightbf}%
  \textfont\bffam=\eightbf \scriptfont\bffam=\sixbf
   \scriptscriptfont\bffam=\fivebf
\normalbaselineskip=12pt
  \setbox\strutbox=\hbox{\vrule height8.5pt depth3.5pt width0pt\relax}%
  \normalbaselines\rm}

%%%%%%%%%%%%%%%%%%%%%%%%%%%%
%% Font families for bold math in title and section heads

\def\eighteenbold{\def\rm{\fam0\eighteenbf}%
  \textfont0=\eighteenbf \scriptfont0=\twelvebf \scriptscriptfont0=\tenbf
  \textfont1=\eighteenmib \scriptfont1=\twelvemib\scriptscriptfont1=\tenmib
  \textfont2=\eighteensyb \scriptfont2=\twelvesyb\scriptscriptfont2=\tensyb
  \textfont3=\eighteenex \scriptfont3=\tenex \scriptscriptfont3=\tenex
  \def\bf{\fam\bffam\eighteenbf}%
  \textfont\bffam=\eighteenbf \scriptfont\bffam=\twelvebf
   \scriptscriptfont\bffam=\tenbf
\normalbaselineskip=20pt
  \setbox\strutbox=\hbox{\vrule height13.5pt depth6.5pt width0pt\relax}%
\everymath {\fam0 }
\everydisplay {\fam0 }
  \normalbaselines\rm}

\def\elevenbold{\def\rm{\fam0\elevenbf}%
  \textfont0=\elevenbf \scriptfont0=\eightbf \scriptscriptfont0=\sixbf
  \textfont1=\elevenmib \scriptfont1=\eightmib \scriptscriptfont1=\sixmib
  \textfont2=\elevensyb \scriptfont2=\eightsyb \scriptscriptfont2=\sixsyb
  \textfont3=\elevenex \scriptfont3=\elevenex \scriptscriptfont3=\elevenex
  \def\bf{\fam\bffam\elevenbf}%
  \textfont\bffam=\elevenbf \scriptfont\bffam=\eightbf
   \scriptscriptfont\bffam=\sixbf
\normalbaselineskip=14pt
  \setbox\strutbox=\hbox{\vrule height10pt depth4pt width0pt\relax}%
\everymath {\fam0 }
\everydisplay {\fam0 }
  \normalbaselines\bf}

%%%%%%%%%%%%%%%%%%%%%%%%%%%%%%%%%%%%%%%%%%%%%%%%%%%%%%%%%
%% 2) Page set up
\hsize=31pc
\vsize=48pc

\parindent=22pt
\parskip=0pt

\widowpenalty=10000
\clubpenalty=10000

\topskip=12pt

\skip\footins=20pt
\dimen\footins=3in % maximum footnote height

\abovedisplayskip=6.95pt plus3.5pt minus 3pt
\belowdisplayskip=\abovedisplayskip

%% Output routine

\voffset=7pt\hoffset= .7in%7pt magstep1

\newif\iftitle%!

\def\amheadline{\iftitle%
\hbox to\hsize{\hss\currannalsline\hss}\else\line{\ifodd\pageno
\hfill\thetitle\hfill\llap{\elevenrm\folio}\else\rlap{\elevenrm\folio}
\hfill\theauthors\hfill\fi}\fi}

\headline={\amheadline}%!!!
\footline={\global\titlefalse}
%\output={\bindingoffset\plainoutput}

%%%%%%%%%%%%%%%%%%%%%%%%%%%%%%%%%%%%%%%
% 3) Title page

 %#1= Volume number, #2=year of publication
\def\annalsline#1#2{\vfill\eject
\ifodd\pageno\else % first page of article on right.
\line{\hfill}
\vfill\eject\fi
\global\titletrue
\def\currannalsline{\eightrm %Annals of Mathematics,%ANNALS
{\eightbf#1} (#2), \thepages}}

\def\titleheadline#1{\def\one{#1}\ifx\one\empty\else
\def\thetitle{{%\frenchspacing%
\let\\ \relax\eightsc\uppercase{#1}}}\fi}

\newif\ifshort

\let\shorttitle\titleheadline

\def\onpages#1#2{\def\thepages{#1--#2}}

\def\thismuchskip[#1]{\vskip#1pt}
\def\ilook{\ifx\next[ \let\go\thismuchskip\else
\let\go\relax\vskip1pt\fi\go}

\def\institution#1{\def\theinstitutions{\vbox{\baselineskip10pt
\def\and{{\eightrm and }}
\def\\{\futurelet\next\ilook}\eightsc #1}}}
\let\institutions\institution

\newwrite\auxfile

\def\startingpage#1{\def\one{#1}\ifx\one\empty\global\pageno=1\else
\global\pageno=#1\fi
\theoremcount=0 \eqcount=0 \sectioncount=0
\openin1 \jobname.aux \ifeof1
\onpages{#1}{???}
\else\closein1 \relax\input \jobname.aux
\onpages{#1}{\lastpage}
\fi\immediate\openout\auxfile=\jobname.aux
}

\def\endarticle{\ifRefsUsed\global\RefsUsedfalse%
\else\vskip21pt\theinstitutions%
\nobreak\vskip8pt
%\vbox{\thereceived\therevised}%
\fi%
\write\auxfile{\string\def\string\lastpage{\the\pageno}}}

\outer\def\bye{\endarticle\par \vfill \supereject \end}

% variation on code from amsspt.sty ==>
\def\document{\let\fontlist@\relax\let\alloclist@\relax
 \elevenpoint}%%% add for annals!!!

% <=== end of code varied from amsppt.sty

\newif\ifacks
\long\def\acknowledgements#1{\def\one{#1}\ifx\one\empty\else
\vskip-\baselineskip%
\global\ackstrue\footnote{\ \unskip}{*#1}\fi}

\def\title#1{\titleheadline{#1}
\vbox to80pt{\vfill
\baselineskip=18pt
\parindent=0pt
\overfullrule=0pt
\hyphenpenalty=10000
\everypar={\hskip\parfillskip\relax}
\hbadness=10000
\def\\ {\vskip1sp}
\eighteenbold#1\vskip1sp}}

\newif\ifauthor

\def\author#1{\vskip11pt
\hbox to\hsize{\hss\tenrm By \tensc#1\ifacks\global\acksfalse*\fi\hss}
\ifshort\else\xdef\theauthors{{\eightsc\uppercase{#1}}}\fi%
\vskip21pt\global\authortrue\everypar={\global\authorfalse\everypar={}}}

\def\twoauthors#1#2{\vskip11pt
\hbox to\hsize{\hss%
\tenrm By \tensc#1 {\tenrm and} #2\ifacks\global\acksfalse*\fi\hss}
\ifshort\else\xdef\theauthors{{\eightsc\uppercase{#1 and #2}}}\fi%
\vskip21pt
\global\authortrue\everypar={\global\authorfalse\everypar={}}}

%%%%%%%%%%%%%%%%%%%%%%%%%%%%%%%%
%% 4) Section heads, counters

\newcount\theoremcount
\newcount\sectioncount
\newcount\eqcount

\newif\ifspecialnumon

\def\eqnumber=#1 {\global\eqcount=#1 \global\advance\eqcount by-1\relax}
\def\sectionnumber=#1 {\global\sectioncount=#1
\global\advance\sectioncount by-1\relax}
\def\proclaimnumber=#1 {\global\theoremcount=#1
\global\advance\theoremcount by-1\relax}

\newif\ifsection
\newif\ifsubsection

\def\elevenboldmath#1{$#1$\egroup}
\def\mathbold{\hbox\bgroup\elevenbold\elevenboldmath}

\def\section#1{\global\theoremcount=0
\global\eqcount=0
\ifauthor\global\authorfalse\else%
\vskip18pt plus 18pt minus 6pt\fi%
{\parindent=0pt
\everypar={\hskip\parfillskip}%            !!! remove
\def\\ {\vskip1sp}\elevenpoint\bf%
\ifspecialnumon\global\specialnumonfalse$\rm\spnum$%
\gdef\sectnum{$\rm\spnum$}%
\else\interlinepenalty=10000%
\global\advance\sectioncount by1\relax\the\sectioncount%
\gdef\sectnum{\the\sectioncount}%
\fi. \hskip6pt#1%                          !!!add }} and stop here
\vrule width0pt depth12pt}
\hskip\parfillskip%\break%!
\global\sectiontrue%
\everypar={\global\sectionfalse\global\interlinepenalty=0\everypar={}}%
\ignorespaces

}

%%%%%%%%%%%%%%%%%%%%%%%%%%%%%%%%
%% 5) Equation Macros

\newif\ifspequation

\let\eqno\leqno %automatic left side equation numbers %%!!!remove l-eqno

\newif\ifineqalignno
\let\saveleqalignno\leqalignno                        %%!!!remove l-eqno
\def\leqalignno{\let\eqnu\Eeqnu\saveleqalignno}

\let\eqalignno\leqalignno

\def\sectandeqnum{%
\ifspecialnumon\global\specialnumonfalse
$\rm\spnum$\gdef\eqnum{{$\rm\spnum$}}\else\global\firstlettertrue
\global\advance\eqcount by1
\ifappend\applett\else\the\sectioncount\fi.%
\the\eqcount
\xdef\eqnum{\ifappend\applett\else\the\sectioncount\fi.\the\eqcount}\fi}

\def\eqnu{\leqno{\hbox{\elevenrm\ifspequation\else(\fi\sectandeqnum
\ifspequation\global\spequationfalse\else)\fi}}}      %!!! l-eqno

\def\Speqnu{\global\setbox\leqnobox=\hbox{\elevenrm
\ifspequation\else%
(\fi\sectandeqnum\ifspequation\global\spequationfalse\else)\fi}}

\def\Eeqnu{\hbox{\elevenrm
\ifspequation\else%
(\fi\sectandeqnum\ifspequation\global\spequationfalse\else)\fi}}

\newif\iffirstletter
\global\firstlettertrue
\def\eqletter#1{\global\specialnumontrue\iffirstletter\global\firstletterfalse
\global\advance\eqcount by1\fi
\gdef\spnum{\the\sectioncount.\the\eqcount#1}\eqnu}

%%% Split math
\newbox\leqnobox
\def\outsideeqnu#1{\global\setbox\leqnobox=\hbox{#1}}

\def\eatone#1{}

%% Vertically centers equation number.
\def\dosplit#1#2{\vskip-.5\abovedisplayskip
\setbox0=\hbox{$\let\eqno\outsideeqnu%
\let\eqnu\Speqnu\let\leqno\outsideeqnu#2$}%
\setbox1\vbox{\noindent\hskip\wd\leqnobox\ifdim\wd\leqnobox>0pt\hskip1em\fi%
$\displaystyle#1\mathstrut$\hskip0pt plus1fill\relax
\vskip1pt
\line{\hfill$\let\eqnu\eatone\let\leqno\eatone%
\displaystyle#2\mathstrut$\ifmathqed~~\qed\fi}}%
\copy1
\ifvoid\leqnobox
\else\dimen0=\ht1 \advance\dimen0 by\dp1
\vskip-\dimen0
\vbox to\dimen0{\vfill
\hbox{\unhbox\leqnobox}
\vfill}
\fi}

\everydisplay{\lookforbreak}

\long\def\lookforbreak #1$${\def\mathone{#1}
\expandafter\testforbreak\mathone\splitmath @}

\def\testforbreak#1\splitmath #2@{\def\mathtwo{#2}\ifx\mathtwo\empty%
#1$$%
\ifmathqed\vskip-\belowdisplayskip
\setbox0=\vbox{\let\eqno\relax\let\eqnu\relax$\displaystyle#1$}%
\vskip-\ht0\vskip-3.5pt\hbox to\hsize{\hfill\qed}
\vskip\ht0\vskip3.5pt\fi
\else$$\vskip-\belowdisplayskip
\vbox{\dosplit{#1}{\let\eqno\eatone
\let\splitmath\relax#2}}%
\nobreak\vskip.5\belowdisplayskip
\noindent\ignorespaces\fi}

%% Proof box to be used when proof ends with equation.

\newif\ifmathqed

%%%%%%%%%%%%%%%%%%%%%%%%%%%%%
%% \mtable, Math table to make binary table easily

%% Use:
% \mtable
% &n_1&n_2&n_3&n_4&n_5&n_6\cr
% \Delta_1&M_3&M_2&0&0&0&0\cr
% \Delta_2&0&0&M_1&M_3&0&0\cr
% \endmtable

\newcount\linenum
\newcount\colnum

%++
\def\spline{\omit&\multispan{\the\colnum}{\hrulefill}\cr}
\def\colcounter{\ifnum\linenum=1\global\advance\colnum by1\fi}

\def\everyline{\noalign{\global\advance\linenum by1\relax}%
\ifnum\linenum=2\spline\fi}

\def\mtable{\bgroup\offinterlineskip
\everycr={\everyline}\global\linenum=0
\halign\bgroup\vrule height 10pt depth 4pt width0pt
\hfill$##$\hfill\hskip6pt\ifnum\linenum>1
\vrule\fi&&\colcounter\hskip12pt\hfill$##$\hfill\hskip12pt\cr}

\def\endmtable{\crcr\egroup\egroup}

%%%%%%%%%%%%%%%%%%%%%%%%%%%%%
% Array

%% Will work in math or in text, will be in math mode inside array.
%% For each column desired supply
%% r, l, or c, for right, left, or center orientation of that column.
%% End each line with \\.

%% To use:
%  \array ccc*
%  x_s\leq a_1\\
%  a_s<x_s^s<b_s\\
%  x_s\geq a_1
%  \endarray

\def\xast{*}
\newcount\intable
\newcount\mathcol
\newcount\savemathcol
\newcount\topmathcol
\newdimen\arrayhspace
\newdimen\arrayvspace

\arrayhspace=8pt % horizontal space between columns, (half this width
                 %  will horizontally precede and follow the array)
\arrayvspace=12pt % vertical space between lines

\newif\ifdollaron

\def\mathalign#1{\def\arg{#1}\ifx\arg\xast%
\let\go\relax\else\let\go\mathalign%
\global\advance\mathcol by1 %
\global\advance\topmathcol by1 %
\expandafter\def\csname  mathcol\the\mathcol\endcsname{#1}%
\fi\go}

\def\arraypickapart#1]#2*{\if#1c \ifmmode\vcenter\else
\global\dollarontrue$\vcenter\fi\else%
\if#1t\vtop\else\if#1b\vbox\fi\fi\fi\bgroup%
\def\one{#2}}

\def\arraystrut{\vrule height .7\arrayvspace depth .3\arrayvspace width 0pt}

\def\array#1#2*{\def\firstarg{#1}%
\if\firstarg[ \def\two{#2} \expandafter\arraypickapart\two*\else%
\ifmmode\vcenter\else\vbox\fi\bgroup \def\one{#1#2}\fi%
\global\everycr={\noalign{\global\mathcol=\savemathcol\relax}}%
\def\\ {\cr}%
\global\advance\intable by1 %
\ifnum\intable=1 \global\mathcol=0 \savemathcol=0 %
\else \global\advance\mathcol by1 \savemathcol=\mathcol\fi%
\expandafter\mathalign\one*%
\mathcol=\savemathcol %
\halign\bgroup&\hskip.5\arrayhspace\arraystrut%
\global\advance\mathcol by1 \relax%
\expandafter\if\csname mathcol\the\mathcol\endcsname r\hfill\else%
\expandafter\if\csname mathcol\the\mathcol\endcsname c\hfill\fi\fi%
$\displaystyle##$%
\expandafter\if\csname mathcol\the\mathcol\endcsname r\else\hfill\fi\relax%
\hskip.5\arrayhspace\cr}

\def\endarray{\crcr\egroup\egroup%
\global\mathcol=\savemathcol %
\global\advance\intable by -1\relax%
\ifnum\intable=0 %
\ifdollaron\global\dollaronfalse $\fi
\loop\ifnum\topmathcol>0 %
\expandafter\def\csname  mathcol\the\topmathcol\endcsname{}%
\global\advance\topmathcol by-1 \repeat%
\global\everycr={}\fi%
}

\def\big#1{{\hbox{$\left#1\vbox to 10pt{}\right.\n@space$}}}
\def\Big#1{{\hbox{$\left#1\vbox to 13pt{}\right.\n@space$}}}
\def\bigg#1{{\hbox{$\left#1\vbox to 16pt{}\right.\n@space$}}}
\def\Bigg#1{{\hbox{$\left#1\vbox to 19pt{}\right.\n@space$}}}

%%%%%%%%%%%%%%%%%%%%%%%%%%%%%%%%%%%%%%%%%%%%%%%%%%%%%%%%%%%%%%%%
% 6) Figure and Table Captions.

\def\figcaption#1#2#3{\topinsert
\vskip4pt %<===topadjust to match height of ascenders on opposing page.
\vbox to#3{\vfill}\vskip1sp
\setbox0=\hbox{\eightsc Figure #1.\hskip12pt\eightpoint #2}
\ifdim\wd0>\hsize
\noindent\eightsc Figure #1.\hskip12pt\eightpoint #2
\else
\centerline{\eightsc Figure #1.\hskip12pt\eightpoint #2}
\fi
\vskip16pt
\endinsert}

\def\wfig#1#2#3{\topinsert
\vskip4pt %<===topadjust to match height of ascenders on opposing page.
\hbox to\hsize{\hss\vbox{\hrule height .25pt width #3
\hbox to #3{\vrule width .25pt height #2\hfill\vrule width .25pt height #2}
\hrule height.25pt}\hss}
\vskip1sp
\centerline{\eightsc Figure #1}
\vskip16pt
\endinsert}

\def\wfigcaption#1#2#3#4{\topinsert
\vskip4pt %<===topadjust to match height of ascenders on opposing page.
\hbox to\hsize{\hss\vbox{\hrule height .25pt width #4
\hbox to #4{\vrule width .25pt height #3\hfill\vrule width .25pt height #3}
\hrule height.25pt}\hss}
\vskip1sp
\setbox0=\hbox{\eightsc Figure #1.\hskip12pt\eightpoint\rm #2}
\ifdim\wd0>\hsize
\noindent\eightsc Figure #1.\hskip12pt\eightpoint\rm #2\else
\centerline{\eightsc Figure #1.\hskip12pt\eightpoint\rm #2}\fi
\vskip16pt
\endinsert}

\def\tabcaption#1#2{\vskip6pt
\setbox0=\hbox{\eightsc Table #1.\hskip12pt\eightpoint #2}
\ifdim\wd0>\hsize
\noindent\eightsc Table #1.\hskip12pt\eightpoint #2
\else
\centerline{\eightsc Table #1.\hskip12pt\eightpoint #2}
\fi
\vskip6pt}

\def\endinsert{\egroup\if@mid\dimen@\ht\z@\advance\dimen@\dp\z@
\advance\dimen@ 12\p@\advance\dimen@\pagetotal\ifdim\dimen@ >\pagegoal
\@midfalse\p@gefalse\fi\fi\if@mid\smallskip\box\z@\bigbreak\else
\insert\topins{\penalty 100 \splittopskip\z@skip\splitmaxdepth\maxdimen
\floatingpenalty\z@\ifp@ge\dimen@\dp\z@\vbox to\vsize {\unvbox \z@
\kern -\dimen@ }\else\box\z@\nobreak\smallskip\fi}\fi\endgroup}

\def\pagecontents{
\ifvoid\topins \else\iftitle\else
\unvbox \topins \fi\fi \dimen@ =\dp \@cclv \unvbox
\@cclv
\ifvoid\topins\else\iftitle\unvbox\topins\fi\fi
\ifvoid \footins \else \vskip \skip \footins \footnoterule
\unvbox \footins \fi \ifr@ggedbottom \kern -\dimen@ \vfil \fi}

%%%%%%%%%%%%%%%%%%%%%%%%%%%%%%%%%%%%%%%%%%%%%%%%%%%%%%%%%%%%%%%%
% 7) End Matter

\newif\ifappend

\def\appendix#1#2{\def\applett{#1}\def\two{#2}%
\global\appendtrue
\global\theoremcount=0
\global\eqcount=0
\vskip18pt plus 18pt
\vbox{\parindent=0pt
\everypar={\hskip\parfillskip}
\def\\ {\vskip1sp}\elevenbold Appendix%
\ifx\applett\empty\gdef\applett{A}\ifx\two\empty\else.\fi%
\else\ #1.\fi\hskip6pt#2\vskip12pt}%
\global\sectiontrue%
\everypar={\global\sectionfalse\everypar={}}\nobreak\ignorespaces}

\newif\ifRefsUsed
\long\def\references{\global\RefsUsedtrue\vskip21pt
\theinstitutions
\global\everypar={}\global\bibnum=0
\vskip20pt\goodbreak\bgroup
\vbox{\centerline{\eightsc References}\vskip6pt}%
\ifdim\maxbibwidth>0pt
\leftskip=\maxbibwidth%
\parindent=-\maxbibwidth%
\else
\leftskip=18pt%
\parindent=-18pt%
\fi
\ninepoint
\frenchspacing
\nobreak\ignorespaces\everypar={\amref}%
}

\def\endreferences{\vskip1sp\egroup\global\everypar={}%
\nobreak\vskip8pt\vbox{\thereceived\therevised}
}

\newcount\bibnum

\def\amref#1 {\global\advance\bibnum by1%
\immediate\write\auxfile{\string\expandafter\string\def\string\csname
\space #1croref\string\endcsname{[\the\bibnum]}}%
\leavevmode\hbox to18pt{\hbox to13.2pt{\hss[\the\bibnum]}\hfill}}

\def\bibline{\hbox to30pt{\hrulefill}\/\/}

\def\name#1{{\eightsc#1}}

\newdimen\maxbibwidth
\def\AuthorRefNames [#1] {%
\immediate\write\auxfile{\string\def\string\cite\string##1{[\string##1]}}

\def\amref{\spamref}
\setbox0=\hbox{[#1] }\global\maxbibwidth=\wd0\relax}

\def\spamref[#1] {\leavevmode\hbox to\maxbibwidth{\hss[#1]\hfill}}

%%%%%%%%%%%%%%%%%%%%%%%%%%%%%%%%%%%%%%%%%%%%%%%%%%%%%%%%%%%%%%%%
%% 8) Footnotes

\def\footnoterule{\kern-3pt\hrule width1in height.5pt\kern2.5pt}

\def\footnote#1#2{%
\plainfootnote{#1}{{\eightpoint\normalbaselineskip11pt
\normalbaselines#2}}}

\def\vfootnote#1{%
\insert \footins \bgroup \eightpoint\baselineskip11pt
\interlinepenalty \interfootnotelinepenalty
\splittopskip \ht \strutbox \splitmaxdepth \dp \strutbox \floatingpenalty
\@MM \leftskip \z@skip \rightskip \z@skip \spaceskip \z@skip
\xspaceskip \z@skip
{#1}$\,$\footstrut \futurelet \next \fo@t}

%%%%%%%%%%%%%%%%%%%%%%%%%%%%%%%%%%%%%%%%%%%%%%%%%%%%%%%%%%%%%%%%
%% 9) Theorem type environments

\newif\iffirstadded
\newif\ifadded

\def\addedlett{}

\def\alltheoremnums{%
\ifspecialnumon\global\specialnumonfalse
\ifadded\global\addedfalse
\iffirstadded\global\firstaddedfalse
\global\advance\theoremcount by1 \fi
\ifappend\applett\else\the\sectioncount\fi.\the\theoremcount\addedlett%
\xdef\theoremnum{\ifappend\applett\else\the\sectioncount\fi.%
\the\theoremcount\addedlett}%
\else$\rm\spnum$\def\theoremnum{{$\rm\spnum$}}\fi%
\else\global\firstaddedtrue
\global\advance\theoremcount by1
\ifappend\applett\else\the\sectioncount\fi.\the\theoremcount%
\xdef\theoremnum{\ifappend\applett\else\the\sectioncount\fi.%
\the\theoremcount}\fi}

\def\allcorolnums{%
\ifspecialnumon\global\specialnumonfalse
\ifadded\global\addedfalse
\iffirstadded\global\firstaddedfalse
\global\advance\corolcount by1 \fi
\the\corolcount\addedlett%
\else$\rm\spnum$\def\corolnum{$\rm\spnum$}\fi%
\else\global\advance\corolcount by1
\the\corolcount\fi}

%% use for Theorem, Corollary, Lemma, Proposition, Demonstration and similar.

\newcount\corolcount
\def\xcorol{Corollary}
\def\xtheorem{Theorem}
\def\xmaintheorem{Main Theorem}

\newif\ifthtitle

\let\saverparen)
\let\savelparen(
\def\rmparenl{{\rm(}}
\def\rmparenr{{\rm\/)}}
{
\catcode`(=13
\catcode`)=13
\gdef\makeparensRM{\catcode`(=13\catcode`)=13\let(=\rmparenl%
\let)=\rmparenr%
\everymath{\let(\savelparen%
\let)\saverparen}%
\everydisplay{\let(\savelparen%
\let)\saverparen\lookforbreak}}}

\medskipamount=8pt plus.1\baselineskip minus.05\baselineskip

\def\rmtext#1{\hbox{\rm#1}}

\def\proclaim#1{\vskip-\lastskip
\def\one{#1}\ifx\one\xtheorem\global\corolcount=0\fi
\ifsection\global\sectionfalse\vskip-6pt\fi
\medskip
{\elevensc#1}%
\ifx\one\xmaintheorem\global\corolcount=0
\gdef\theoremnum{Main Theorem}\else%
\ifx\one\xcorol\ \allcorolnums\else\ \alltheoremnums\fi\fi%
\ifthtitle\ \global\thtitlefalse{\rm(\thethtitle)}\fi.%
\hskip1em\bgroup\let\text\rmtext\makeparensRM\it\ignorespaces}

\def\nonumproclaim#1{\vskip-\lastskip
\def\one{#1}\ifx\one\xtheorem\global\corolcount=0\fi
\ifsection\global\sectionfalse\vskip-6pt\fi
\medskip
{\elevensc#1}.\ifx\one\xmaintheorem\global\corolcount=0
\gdef\theoremnum{Main Theorem}\fi\hskip.5pc%
\bgroup\it\makeparensRM\ignorespaces}

\def\endproclaim{\egroup\medskip}

%% Use demo for Proof, Proof of, Definition, Example,
%% Remark, Case, Subcase, Conjecture, Note, Notation,
%% Convention, Construction and Step.
%% Any other use for demo will format similar to `Proof.'

\def\xproof{Proof}
\def\xremark{Remark}
\def\xcase{Case}
\def\xsubcase{Subcase}
\def\xconjecture{Conjecture}
\def\xstep{Step}
\def\xof{of}

\def\deconstruct#1 #2 #3 #4 #5 @{\def\one{#1}\def\two{#2}\def\three{#3}%
\def\four{#4}%
\ifx\two\empty #1\else%
\ifx\one\xproof%
\ifx\two\xof%
  \ifx\three\xcorol Proof of Corollary \rm#4\else%
     \ifx\three\xtheorem Proof of Theorem \rm#4\else\xone\fi%
  \fi\fi%
\else\xone\fi\fi.}

\def\pickup#1 {\def\this{#1}%
\ifx\this\xproof\global\let\go\demoproof
\global\let\enddemo\endproof\else
\ifx\this\xremark\global\let\go\demoremark\else
\ifx\this\xcase\global\let\go\demostep\else
\ifx\this\xsubcase\global\let\go\demostep\else
\ifx\this\xconjecture\global\let\go\demostep\else
\ifx\this\xstep\global\let\go\demostep\else
\global\let\go\demoproof\fi\fi\fi\fi\fi\fi}

\newif\ifnonum
\def\demo#1{\vskip-\lastskip
\ifsection\global\sectionfalse\vskip-6pt\fi
\def\one{#1 }\def\two{#1*}%
\setbox0=\hbox{\expandafter\pickup\one}\expandafter\go\two}

\def\numbereddemo#1{\vskip-\lastskip
\ifsection\global\sectionfalse\vskip-6pt\fi
\def\two{#1*}%
\expandafter\demoremark\two}

\def\demoproof#1*{\medskip\def\xone{#1}
{\ignorespaces\it\expandafter\deconstruct\xone {} {} {} {} {} @%
\unskip\hskip6pt}\rm\ignorespaces}

\def\demoremark#1*{\medskip
{\it\ignorespaces#1\/} \ifnonum\global\nonumtrue\else
 \alltheoremnums\unskip.\fi\hskip1pc\rm\ignorespaces}

\def\demostep#1 #2*{\vskip4pt
{\it\ignorespaces#1\/} #2.\hskip1pc\rm\ignorespaces}

\def\enddemo{\medskip}

\def\endproof{\ifmathqed\global\mathqedfalse\medskip\else
\parfillskip=0pt~~\hfill\qed\medskip
\fi\global\parfillskip0pt plus 1fil\relax
\gdef\enddemo{\medskip}}

\def\qed{\vbox{\hrule\hbox{\vrule height6pt\hskip6pt\vrule}\hrule}}

%% Proof box to be used in a \proclaim{}...\endproclaim environment

%%%%

%%%%%%%%%%%%%%%%%%%%%
%% 10) CrossRefs

%%% Generic crossreferencing
%%% to use: \label\nameoflabel* (will give the page number when referenced)

% Commands to access current state of counter, for cross-referencing
% \sectnum
% \theoremnum
% \eqnum

%%% You can make another definition that includes counters and/or the
%%% page number and access this information as the second argument:
%%% \label\yourlabelname[2.13]*

%%% Since this method of cross-referencing relies
%%% on an auxiliary file, the first time you tex the file
%%% you will get `??' when you write \ref\nameoflabel.
%%% When you TeX the file the second time the auxiliary file
%%% will be input and \ref\nameoflabel will produce the cross-ref.

\def\stripbs#1#2*{\def\one{#2}}

\def\emptyspace{ }
\def\nextthing{}
\def\newline{***}
\def\eatone#1{ }

\def\lookatspace#1{\ifcat\noexpand#1\ \else%
\gdef\nextthing{}\xdef\next{#1}%
\ifx\next\emptyspace%
\let\nextthing\emptyspace\else\ifx\next\newline%
\gdef\nextthing{\eatone}\fi\fi\fi\egroup\nextthing#1}

{\catcode`\^^M=\active%
\gdef\spacer{\bgroup\catcode`\^^M=\active%
\let^^M=\newline\obeyspaces\lookatspace}}

\def\ref#1{\seeifdefined{#1}\expandafter\csname\one\endcsname\spacer}

\def\cite#1{\expandafter\ifx\csname#1croref\endcsname\relax[??]\else
\csname#1croref\endcsname\fi\spacer}

%% for testing in \label and \ref to see if term already labeled.

\def\seeifdefined#1{\expandafter\stripbs\string#1croref*%
\crorefdefining{#1}}

\newif\ifcromessage
\global\cromessagetrue

\def\crorefdefining#1{\ifdefined{\one}{}
{\ifcromessage\global\cromessagefalse%
\message{\spaces\spaces\spaces\spaces\spaces\spaces\spaces}%
\message{<Undefined reference.}%
\message{Please TeX file once more to have accurate cross-references.>}%
\message{\spaces\spaces\spaces\spaces\spaces\spaces\spaces}\fi[??]}}

\def\label#1#2*{\gdef\ctest{#2}%
\xdef\currlabel{\string#1croref}
\expandafter\seeifdefined{#1}%
\ifx\empty\ctest%
\xdef\labelnow{\write\auxfile{\noexpand\def\currlabel{\the\pageno}}}%
\else\xdef\labelnow{\write\auxfile{\noexpand\def\currlabel{#2}}}\fi%
\labelnow}

\def\ifdefined#1#2#3{\expandafter\ifx\csname#1\endcsname\relax%
#3\else#2\fi}

%%%%%%%%%%%%%%%%%%%%%
%% 11) Listing

%% To use with asterisks:

%%%%%%%%%%%%%%%%%%%%%%
%% 12) Article and Journal Table of Contents

\def\articlecontents{
\vskip20pt\centerline{\bf Table of Contents}\everypar={}\vskip6pt
\bgroup \leftskip=3pc \parindent=-2pc
\def\item##1{\vskip1sp\indent\hbox to2pc{##1.\hfill}}}

\def\endcontents{\vskip1sp\leftskip=0pt\egroup}

\def\journalcontents{\vfill\eject
\def\currannalsline{\hfill}
\global\titletrue
\vglue3.5pc
\centerline{\tensc\hskip12pt TABLE OF CONTENTS}\everypar={}\vskip30pt
\bgroup \leftskip=34pt \rightskip=-12pt \parindent=-22pt
  \def\\ {\vskip1sp\noindent}
\def\pagenum##1{\unskip\parfillskip=0pt\dotfill##1\vskip1sp
\parfillskip=0pt plus 1fil\relax}
\def\name##1{{\tensc##1}}}

%% default values

\institution{}
\onpages{0}{0}
\def\lastpage{???}
\def\thetitle{Title ???}
\def\theauthors{Authors ???}
\def\thereceived{}
\def\therevised{}

\gdef\split{\relaxnext@\ifinany@\let\next\insplit@\else
 \ifmmode\ifinner\def\next{\onlydmatherr@\split}\else
 \let\next\outsplit@\fi\else
 \def\next{\onlydmatherr@\split}\fi\fi\let\eqnu\xspliteqnu\next}

\gdef\align{\relaxnext@\ifingather@\let\next\galign@\else
 \ifmmode\ifinner\def\next{\onlydmatherr@\align}\else
 \let\next\align@\fi\else
 \def\next{\onlydmatherr@\align}\fi\fi\let\eqnu\xspliteqnu\next}

\def\spliteqnu{{\tenrm\sectandeqnum}\relax}

\def\xspliteqnu{\tag\spliteqnu}

\catcode`@=12

\document
\sectionfalse
%-------------- Publisher's entries --------------------
\annalsline{November}{1995}
%\line{\hfil Revised }
\startingpage{1}     %numeration
%%\received{??}
%\revised{February, 1995}

\comment
\nopagenumbers
\headline{\ifnum\pageno=1\hfil\else \rightheadline\fi}
%{\ifodd\pageno\rightheadline \else \leftheadline\fi}\fi}%%after else
\def\rightheadline{\hfil\eightit
%! Running title (odd page)
Kadell conjectures
\quad\eightrm\folio}

\voffset=2\baselineskip
\endcomment

%\magnification=\magstep1

%--------------- Author macros ---------------
%                   MACROS
%
%                                 AUX
%
%
%
%                      endaux
%

\def\for{\  \hbox{ for } \ }
\def\if{ \ \hbox{ if } \ }

\def\where{\  \hbox{ where } \ }
\def\and{\  \hbox{ and } \ }

\def\equal{\buildrel  def \over =}

\def\om{\omega}

\def\th{\theta}
\def\al{\alpha}

\def\de{\delta}
\def\De{\Delta}

\def\Ga{\Gamma}
\def\ze{\zeta}

    %from copy, ell

\def\vep{\varepsilon}

\def\tal{\tilde{\alpha}}

\def\tw{\tilde w}

\def\tz{\tilde z}
\def\tb{\tilde b}

\def\hT{\hat{T}}

\def\hw{\hat{w}}

\def\hv{\hat{v}}

\def\C{\bold{C}}
\def\Q{\bold{Q}}

\def\R{\bold{R}}

\def\Z{\bold{Z}}

\def\one{\bold{1}}

\def\0{\bold{0}}

\def\C{\hbox{\bf C}}

%\def\bs{\hbox{\bf S}}          %ell
% macdonald

\def\l{\Cal{L}}

\def\w{\Cal{W}}

\font\germ=eufb10 %at 12pt
%\font\germm=eufb9 at 12pt
%\font\germ=eufm9 at 12pt
\def\goth#1{\hbox{\germ #1}}

\def\TT{\goth{T}}

\def\AA{\goth{A}}

\font\smm=msbm10 at 12pt
\def\symbol#1{\hbox{\smm #1}}
\def\lsmash{{\symbol n}}

%endmacros

%------------------------------------------------------------------
%-------------- Author entries --------------------

%\comment                               %to remove the title
\title
{Kadell's two conjectures \\ for Macdonald polynomials}

 %Article title
\shorttitle{ Kadell Conjectures}
 % Shortened version for headline title

% Acknowledgements: Please enter all acknowledgements here.
\acknowledgements{
Partially supported by NSF grant DMS--9301114}

% Please uncomment and use appropriate command:
\author{ Ivan Cherednik}
%\twoauthors{}{}
%\authors{}% Separate each author with a comma and a space.

% Institution:
\institutions{
Math. Dept, University of North Carolina at Chapel Hill,
 N.C. 27599-3250
\\ Internet: chered\@math.unc.edu
}

%\endcomment                              %to remove the title
%-------------- Article Text--------------------

%\intro %(Optional, Introduction)
%
%
%
%                        INTRO
%
%
%{\bf 0. Introduction.}
\vfil
%\sectionfalse

Recently  Kevin Kadell found  interesting properties
of anti-symmetric variants of the so-called Jack polynomials [Ka].
He formulated two conjectures about negative integral and
half-integral values of the parameter $k$ ($k=1$ for
the characters of compact simple Lie groups).
As it was observed independently by Ian Macdonald
and the author,
 these conjectures follow readily
from the interpretation of the Jack polynomials as eigenfunctions
of the Calogero- Sutherland -Heckman -Opdam
 second order operators  generalizing the radial
parts of the Laplace operators on symmetric spaces (see
[HO,He,M1,M2]).
  The difference case requires a bit different
treatment but still is not complicated. We will formulate and
prove the Kadell conjectures for the Macdonald polynomials
(the $q,t$-case).

These statements are of  certain interest
because negative $k$ are somehow connected with irreducible
represenations for anti-dominant highest weights (and with
represenations of Kac- Moody algebras of negative integral
central charge). They also make more complete the theory
of Macdonlad's polynomials at roots of unity started in
[Ki], [C3,C4].  Half-integral $k$ appear in the theory of spherical
functions. Anyway it is challenging to understand what is
going on when $0> k\in \Q$, since  these values are
singular for the coefficients of symmetric Macdonald polynomials.

The author thanks Kevin Kadell for explaining his conjectures
and Ian Macdonald who stimulated this note a lot.
It was started in San Diego. The author is grateful to Adriano
Garsia for the kind invitation and hospitality.

\section{Jack polynomials}
Let $R=\{\al\}   \subset \R^n$ be a root system of type $A,B,...,F,G$
with respect to a euclidean form $(z,z')$ on $\R^n \ni z,z'$,
$W$ the Weyl group  generated by the the reflections $s_\al$.
We assume that $(\al,\al)=2$ for long $\al$.
Let us  fix the set $R_{+}$ of positive  roots ($R_-=-R_+$),
the corresponding simple
roots $\al_1,...,\al_n$, and  their dual counterparts
$a_1 ,..., a_n,  a_i =\al_i^\vee, \where \al^\vee =2\al/(\al,\al)$.
The dual fundamental weights
$b_1,...,b_n$  are determined from the relations  $ (b_i,\al_j)=
\de_i^j $ for the
Kronecker delta. We will also introduce the dual root system
$R^\vee =\{\al^\vee, \al\in R\}, R^\vee_+$, and the lattices
$$
\eqalignno{
& A=\oplus^n_{i=1}\Z a_i \subset B=\oplus^n_{i=1}\Z b_i,
}
$$
  $A_+, B_+$  for $\Z_{+}=\{m\in\Z, m\ge 0\}$
instead of $\Z$. (In the standard notations, $A= Q^\vee,\
B = P^\vee $ - see [B].)  Later on,
$$
\eqalign{
&\nu_{\al}=\nu_{\al^\vee}\ =\ (\al,\al),\  \nu_i\ =\ \nu_{\al_i}, \
\nu_R\ = \{\nu_{\al}, \al\in R\}, \cr
&\rho_\nu\ =\ (1/2)\sum_{\nu_{\al}=\nu} \al \ =
\ (\nu/2)\sum_{\nu_i=\nu}  b_i, \for\al\in R_+,\cr
&r_\nu\ =\ \rho_\nu^\vee \ =\ (2/\nu)\rho_\nu\ =\
\sum_{\nu_i=\nu}  b_i,\quad 2/\nu=1,2,3.
}
\eqnu
$$

We set $x_i=exp({b_i}),\  x_b=exp(b)= \prod_{i=1}^n
x_i^{\kappa_i} \for b=\sum_{i=1}^n \kappa_i b_i$,
 $\C[x]\equal\C[x_b, b\in B]$,  $\partial_a(x_b)=(a,b)x_b$.
The
 monomial symmetric functions
$m_{b}\ =\ \sum_{c\in W(b)}x_{-c}$ for $b\in B_+$
form a base of the space
 $\C[x]^W$ of all $W$-invariant polynomials (note
the sign of $c$).

Let $ k_\nu\in \C, \nu\in \nu_R$,
$k_\al=k_{\nu_\al}=k_{\al^{\vee}}, r_k=\sum_{\nu}k_\nu r_\nu$.
The operator
$$
\eqalignno{
&L^{(k)}_2\ =\ \sum_{i=1}^n\partial_{\al_i}\partial_{b_i}
+\sum_{a\in R_+^\vee}  k_a {x_a+1\over x_a-1}\partial_a
+(r_k,r_k),
&\eqnu
\label\ljack\eqnum*
}
$$
 can be also represented as $L^{(k)}_2=\De_k^{-1}H^{(k)}_2\De_k$
for
$$
\eqalignno{
H^{(k)}_2\ =\ \sum_{i=1}^n\partial_{\al_i}\partial_{b_i}
+&\sum_{a\in R_+^\vee}  k_a(1-k_a)(a,a)
(\exp(a/2)-\exp(-a/2))^{-2},\cr
&\De_k\ =\ \prod_{a\in R_+^\vee}(\exp(a/2)-\exp(-a/2))^{k_a}.
 &\eqnu
\label\hjack\eqnum*
}
$$
Here $\De_k$ for non-integral $k$ should
be a
 solution of the
obvious system of differential equations. The passage from
$L$ to $H$ is well known (see e.g. [C1], Corollary 2.8).

Jack polynomials $j^{(k)}_b(x)\in \C[x], b\in B_+$ belong to the
algebra $\C[x]^W$ of W-invariant polynomials. They can be
fixed up to proportionality from the following eigenvalue problem:
$$
\eqalignno{
&L^{(k)}_2 (j^{(k)}_b)\ =\ (b+r_k,b+r_k) j^{(k)}_b.
&\eqnu
\label\ejack\eqnum*}
$$
More exactly (see e.g. [He, M2, O]),
they are determined uniquely by means of this equation and
 the conditions
$$
\eqalignno{
&J^{(k)}_b-m_b\ \in\ \oplus_c\C(q,t)m_{c} \for c\prec b  &\eqnu
 \cr
&\hbox{where\ }
 c\in B_+,\
 c\prec b \hbox{\ means\ that\ }  b-c \in A_+, c\neq b.
\label\jack\eqnum*
}
$$
Here $k_\nu$ are arbitrary except negative rational
numbers.

A more traditional approach  is as follows
(see also [M1] and the papers by Hanlon, Stanley).
Let $\langle f\rangle$ be the constant term of $f\in \C[x]$.
Setting for $k_\nu \in \Z$,
$$
\eqalignno{
&\langle f,g\rangle\ =\langle \De_{2k} f(x)\ g(x^{-1})\rangle\ =\
\langle g,f\rangle \where
f,g \in \C(q,t)[x]^W,
&\eqnu
\label\jackpro\eqnum*
}
$$
one can  introduce  $J^{(k)}_b$   by means of
the conditions $\langle j^{(k)}_b, m_{c}\rangle = 0, \for c\prec b $
together with
(\ref\jack).
They are pairwise orthogonal for arbitrary $b\in B_+$,
since $H^{(k)}_2$ is self-adjoint and therefore
$L^{(k)}_2$ is self-adjoint with respect to the pairing
$\langle\ ,\ \rangle$ (the eigenvalues distinguish
$b$ at least for non-negative $k$). If $k$ is not an integer then the
pairing should be understood analytically.

\proclaim {Theorem (Kadell Conjectures)}
Let  $m_\nu\in \Z_+$. Then
$$
\De_{2m} j^{(1/2 + m)}_b\ =\ j^{(1/2 - m)}_{b+2r_m},\
b\in B_+.
$$
The polynomial
$\De_{2m+1} j^{(m+1)}_b$ is  anti-symmetric
 ($s_i(\  )= -(\ )$, $1\le i\le n$)
and satisfies  (\ref\ejack)
for $k=-m$ and the eigenvalue equal to $(b+r_{m+1},b+r_{m+1})$.
\endproclaim
The proof results immediately from the invariance of $H^{(k)}_2$
when $k$ is replaced by $1-k$.

\section { Macdonald polynomials}
Let us
introduce the algebra
$\C(q,t)[x]$  of polynomials in terms of $x_i^{\pm 1}$ with the
coefficients belonging to the field $\C(q,t)$ of rational functions
in terms  of indefinite complex parameters $q, t_\nu,
\nu\in \nu_R$ (we will put $t_\al=t_{\nu_\al}=t_{\al^\vee}$).
The coefficient of $x^0=1$ ({ the constant term})
will be again denoted by $\langle \  \rangle$. The following product is a
Laurent series in $x$ with the coefficients in  the algebra
$\C[t][[q]]$ of formal series in $q$ over polynomials in $t$:
$$
\eqalign{
&\mu=\mu_{q,t}\ =\ \prod_{a \in R_+^\vee}
\prod_{i=0}^\infty {(1-x_aq_a^{i}) (1-x_a^{-1}q_a^{i+1})
\over
(1-x_a t_a q_a^{i}) (1-x_a^{-1}t_a^{}q_a^{i+1})},
}
\eqnu
\label\mu\eqnum*
$$
where $q_a=q_{\nu}=q^{2/\nu} \for \nu=\nu_a$.
We note that  $\mu\in
\C(q,t)[x]$ if $t_\nu=q_\nu^{k_\nu}$ for $k_\nu\in \Z_+$.
The coefficients of $\mu_1\equal \mu/\langle \mu \rangle$
are from $\C(q,t)$, where the formula for the
constant term of $\mu$ is as follows
(see [C2]):
$$
\eqalign{
&\langle\mu\rangle\ =\ \prod_{a \in R_+^\vee}
\prod_{i=1}^\infty {(1-x_a(t^\rho)q_a^{i})^2
\over
(1-x_a(t^\rho) t_aq_a^{i}) (1-x_a(t^\rho) t_a^{-1}q_a^{i})}.
}
\eqnu
\label\consterm\eqnum*
$$
%Here and further
%$x_{b}(t^{\pm\rho}q^{c})=
%q^{(b,c)}\prod_\nu t_\nu^{\pm(b,\rho_\nu)}$.
We note that
$\mu_1^*\ =\ \mu_1$  with respect to the involution
$$
 x_b^*\ =\  x_{-b},\ t^*\ =\ t^{-1},\ q^*\ =\ q^{-1}.
$$

Setting ,
$$
\eqalignno{
&\langle f,g\rangle\ =\langle \mu_1 f\ {g}^*\rangle\ =\
\langle g,f\rangle^* \for
f,g \in \C(q,t)[x]^W,
&\eqnu
\label\innerpro\eqnum*
}
$$
 we  introduce the  Macdonald
polynomials $p_b^{q,t}= p_b(x),\   b \in B_+$, by means of
the conditions
$$
\eqalignno{
&p_b-m_b\ \in\ \oplus_c\C(q,t)m_{c},\
\langle p_b, m_{c}\rangle = 0, \for c\prec b.  &\eqnu
\label\macd\eqnum*
}
$$
They can be determined by the Gram - Schmidt process
because the (skew Macdonald) pairing (see [M1,M2,C2])
is non-degenerate (for generic $q,t$)
 and form a
basis in $\C(q,t)[x]^W$. As it was established by Macdonald
they are pairwise orthogonal for arbitrary $b\in B_+$.
We note that $p_b$ are ''real'' with respect to the formal
 conjugation sending $q\to q^{-1},\ t\to t^{-1}$.
It makes our definition compatible with  Macdonald's
original one (his $\mu$ is somewhat different).

The construction is applicable when
 $t_\nu = q_\nu^{k_\nu}$ for $k_\nu\in \Z_+$. More exactly,
the formulas for the coefficients of Macdonald polynomials
have singularities only for rational negative $k$
(if $q$ is generic). We come to
the
 operator interpretation
of the Macdonald polynomials.

\section { Affine root systems}
The vectors $\ \tal=[\al,k] \in
\R^n\times \R \subset \R^{n+1}$
for $\al \in R, k \in \Z $
form the { affine root system}
$R^a \supset R$ ( $z\in \R^n$ are identified with $ [z,0]$).
We add  $\al_0 \equal [-\th,1]$ to the  simple roots
for the { maximal root} $\th \in R$.
The corresponding set $R^a_+$ of positive roots coincides
with $R_+\cup \{[\al,k],\  \al\in R, \  k > 0\}$. See [B].

We denote the Dynkin diagram and its affine completion with
$\{\al_j,0 \le j \le n\}$ as the vertices by $\Ga$ and $\Ga^a$.
Let $m_{ij}=2,3,4,6$\  if $\al_i\and\al_j$ are joined by 0,1,2,3 laces
respectively.
The set of
the indices of the images of $\al_0$ by all
the automorphisms of $\Ga^a$ will be denoted by $O$ ($O=\{0\}
\for E_8,F_4,G_2$). Let $O^*={r\in O, r\neq 0}$.

Given $\tal=[\al,k]\in R^a,  \ b \in B$, let
$$
\eqalignno{
&s_{\tal}(\tz)\ =\  \tz-(z,\al^\vee)\tal,\
\ b'(\tz)\ =\ [z,\ze-(z,b)]
&\eqnu
%&(1.1)
}
$$
for $\tz=[z,\ze] \in \R^{n+1}$.

The { affine Weyl group} $W^a$ is generated by all $s_{\tal}$.
 One can take
the simple reflections $s_j=s_{\al_j}, 0 \le j \le n,$ as its
generators and introduce the corresponding notion of the
length. This group is
the semi-direct product $W\lsmash A'$ of
its subgroups $W$ and $A'=\{a', a\in A\}$, where
$$
\eqalignno{
& a'=\ s_{\al}s_{[\al,1]}=\ s_{[-\al,1]}s_{\al}\for a=\al^{\vee},
\ \al\in R.
&\eqnu
%&(1.2)
}
$$

The { extended Weyl group} $ W^b$ generated by $W\and B'$
(instead of $A'$) is isomorphic to $W\lsmash B'$:
$$
\eqalignno{
&(wb')([z,\ze])\ =\ [w(z),\ze-(z,b)] \for w\in W, b\in B.
&\eqnu
}
$$

 Given $b\in B_+$, let
$$
\eqalignno{
&\om_{b} = w_0w^+_0  \in  W,\ \pi_{b} =
b'(\om_{b_+})^{-1}
\ \in \ W^b,\ \pi_i=\pi_{b_i},
&\eqnu
\label\w0\eqnum*
}
$$
where $w_0$ (respectively, $w^+_0$) is the longest element in $W$
(respectively, in $ W_{b}$ generated by $s_i$ preserving $b$)
relative to the
set of generators $\{s_i\}$ for $i >0$.

We will  use here only the
elements $\pi_r=\pi_{b_r}, r \in O$. They leave $\Ga^a$ invariant
and form a group denoted by $\Pi$,
 which is isomorphic to $B/A$ by the natural
projection $\{b_r \to \pi_r\}$.
The relations $\pi_r(\al_0)= \al_r
$ separate the
indices $r \in O^*$ (see e.g. [C2]), and
$$
\eqalignno{
& W^b  = \Pi \lsmash W^a, \where
  \pi_rs_i\pi_r^{-1}  =  s_j \if \pi_r(\al_i)=\al_j,\  0\le j\le n.
&\eqnu
%&(1.6)
}
$$

Given $\nu\in\nu_R,\  r\in O^*,\  \tw \in W^a$, and a reduced
decomposition $\tw\ =\ s_{j_l}...s_{j_2} s_{j_1} $ with respect to
$\{s_j, 0\le j\le n\}$, we call $l\ =\ l(\hw)$ the { length} of
$\hw = \pi_r\tw \in W^b$.

Later on  $b$ and $b'$ will not be distinguished.
We set $([a,k],[b,l])=(a,b)$ for $a,b\in B,\
 a_0=\al_0,\ \nu_{\al^\vee}=\nu_\al.$

\section { Difference operators}
We put
$m=2 \for D_{2k} \and C_{2k+1},\ m=1 \for C_{2k}, B_{k}$,
otherwise $m=|\Pi|$. Later on $ \C_{ q}$
 is the field of rational
functions in $ q^{1/m}$.
 Setting
$$
\eqalignno{
& x_{\tb}=  \prod_{i=1}^nx_i^{k_i} q^{ k} \if
\tb=[b,k],
b=\sum_{i=1}^nk_i b_i\in B,\ k \in {1\over m}\Z,
&\eqnu
\label\xde\eqnum*}
$$
 we will identify polynomials with the corresponding
  multiplication  operators  in $\C_q[x]=
\C_q[x_1^{\pm 1},\ldots,x_n^{\pm 1}]$.
We replace $ \C_{ q}$
by $ \C_{ q,t}$ if the functions (coefficients)
also depend rationally
on $\{t_\nu^{1/2} \}$.

The elements $\hw \in W^b$ act in $\C_{ q}[x]$
 by the
formulas:
$$
\eqalignno{
&\hw(x_{\tb})\ =\ x_{\hw(\tb)}.
&\eqnu}
$$
 In particular:
$$
\eqalignno{
&\pi_r(x_{b})\ =\ x_{\om^{-1}_r(b)} q^{(b_{r^*},b)}
\for \al_{r^*}\ =\ \pi_r^{-1}(\al_0), \ r\in O^*.
&\eqnu}
$$
\label\pi\eqnum*

The { Demazure-Lusztig operators} (see
[C1,C2] for more detail )
$$
\eqalignno{
&T_j^{q,t} =  t_j ^{1/2} s_j\ +\
(t_j^{1/2}-t_j^{-1/2})(x_{a_j}-1)^{-1}(s_j-1),
\ 0\le j\le n.
&\eqnu
\label\Demaz\eqnum*
}
$$
act   in $\C_{ q,t}[x]$ naturally.
We note that only $\hT_0$ depends on $ q$:
$$
\eqalign{
&T_0\  =  t_0^{1/2}s_0\ +\ (t_0^{1/2}-t_0^{-1/2})
( q x_{\th}^{-1} -1)^{-1}(s_0-1),\cr
&\where
s_0(x_i)\ =\ x_ix_{\th}^{-(b_i,\th)} q^{(b_i,\th)}.
}
%\eqno(2.12)
\eqnu
$$

Given $\tw \in W^a, r\in O,\ $ the product
$$
\eqalignno{
&T_{\pi_r\tw}\equal \pi_r\prod_{k=1}^l T_{i_k},\where
\tw=\prod_{k=1}^l s_{i_k},
l=l(\tw),
&\eqnu
\label\Tw\eqnum*
}
$$
does not depend on the choice of the reduced decomposition
(because $\{T\}$ satisfy the same ``braid'' relations as $\{s\}$ do).
Moreover,
$$
\eqalignno{
&T_{\hv}T_{\hw}\ =\ T_{\hv\hw}\  \hbox{ whenever}\
 l(\hv\hw)=l(\hv)+l(\hw) \for
\hv,\hw \in W^b.
 %&(2.7)}
&\eqnu}
$$
\label\TT\eqnum*
  In particular, we arrive at the pairwise
commutative elements
$$
\eqalignno{
& Y_{b}\ =\  \prod_{i=1}^nY_i^{k_i} \if
b=\sum_{i=1}^nk_ib_i\in B,\where
 Y_i\equal T_{b_i},
&\eqnu
\label\Yb\eqnum*
}
$$
satisfying the relations
$$
\eqalign{
&T^{-1}_iY_b T^{-1}_i\ =\ Y_b Y_{a_i}^{-1} \if (b,\al_i)=1,
\cr
& T_iY_b\ =\ Y_b T_i \if (b,\al_i)=0, \ 1 \le i\le  n.
}
%\eqno(2.9)
\eqnu
$$

Given an operator
$$
\eqalignno{
& H\ = \sum_{b\in B, w\in W} h_{b,w} b  w,
&\eqnu
\label\hatH\eqnum*
}
$$
where  $h_{b,w}$ belong to  the field $\C_{ q,t}(x)$
of rational  functions in
$x_1,...,x_n$, we set
$$
\eqalign{
&[H]_{\dagger}= \sum h_{b,w}  b.
}
\eqnu
\label\Brack\eqnum*
$$
We will use the following  theorem from [C2,C3] (see also
[M3]).

\proclaim{Theorem}
(i) The difference operators $\{ L^{q,t}_f\equal
[f(Y_1,\ldots,Y_n)]_{\dagger}$
for $f\in \C[x]^W\}$
are pairwise commutative,  $W$-invariant (i.e $w L_f w^{-1}=$
$L_f$ for all $w\in W$) and preserve $\C_{ q,t}[x]^W$. The
Macdonald polynomials $p_b=p_b^{ q,t} (b\in B_+)$
from  (\ref\macd) are their eigenvectors:
$$
\eqalignno{
&L_f(p_b^{ q,t})=f(t^\rho q^{b}) p_b^{ q,t},\
x_i(t^\rho q^{b})\equal
 q^{(b_i,b)}\prod_\nu t_\nu^{(b_i,\rho_\nu)}.
&\eqnu
\label\Lf\eqnum*}
$$

(ii) Given  $v\in \nu_R $ such that $\nu\in v$
if $t_\nu\neq 1$, we set
$$
 d_v^{q,t} =  \prod_{\nu_a\in v}((t_a x_{a})^{1/2}-
(t_a x_{a})^{-1/2}),\ tq_v=\{t_\nu \hbox{\ if\ }\nu\not\in v,
\ t_\nu q^{2/\nu}
\hbox{\ otherwise}\}.
$$
The polynomilals $g_b=g^{q,t,v}_b\ =\ d_v^{q,t}p^{q,tq_v}_b$
satisfy the relations
$$
\eqalignno{
&f(Y^{q,t})(g_b)\ =\ f((tq_v)^\rho q^{b}) g_b \for
f\in \C[x]^W,\cr
&T_i^{t}(g_b)\ =\ \vep t_i^{\vep 1/2} g_b,\ 1\le i\le n,
 &\eqnu
\label\Lgv\eqnum*
}
$$
where $\vep=\pm$ if $\nu_i\not\in v,\ \nu_i\in v$
respectively.
\endproclaim
\label\LF\theoremnum*

\comment
\proclaim{Corollary}
If $q$ is generic then the Macdonald polynomials
are well defined when $t_\nu \neq q^{-k_\nu}$ for
all $\nu$ and $k_\nu=1,2,\ldots.
\endproclaim
{\it Proof.}
Let $t_\nu=q^{k_\nu}$ and $k$ are non-negative.
Since $f$ are arbitrary symmetric,
the polynomilal $p_b, b\in B_+$ exists if
$ B_+\ni w(r_k+b)-r_k\prec b$ for all $w\in W$.
Assuming that the latter does not hold,
$(w(r_k)-r_k, \al)\in \Z$ for all $\al\in R$.
Hence
\endcomment

\section { Kadell's formulas}
The relations
$$
\eqalignno{
&(s_i)^\iota = -s_i,\ x_i^\iota=x_i,\  0\le i\le n,\
 q^\iota=q,\ (t^{1/2})^\iota = -t^{-1/2}
 &\eqnu
\label\iota\eqnum*
}
$$
can be extended to an involution on functions of $x$
and operators. We will also use the conjugation
which sends $t^{1/2}$ to $-t^{1/2}$ (and does not change other
generators of $\C_{q,t}[x]$). The latter fixes the Macdonald
polynomials since they are expressed via $t$ only.
We write $p^{(k)}, L^{(k)}$ and so on if $t_\nu=q_\nu^{k_\nu}$.

\proclaim{ Main Theorem}
(i) Let  $m_\nu\in \Z_+, \ t_\nu\ =\ q_\nu^{m_\nu+1/2}$
($\nu\in \nu_R$). Then
$$
\eqalign{
&\de_{2m} p^{(1/2 + m)}_b\ =\ p^{(1/2 - m)}_{b+2r_m} \for
\cr
&\de_{2m}\  =\ \prod_{a \in R_+^\vee}
 \{(x_a^{1/2} t_a^{1/2}q_a^{-1/2}- x_a^{-1/2} t_a^{-1/2}
q_a^{1/2})\cdots
\cr
&(x_a^{1/2} t_a^{1/2}q_a^{-i/2}- x_a^{-1/2} t_a^{-1/2}q_a^{i/2})
\cdots (x_a^{1/2} t_a^{-1/2}q_a^{1/2}- x_a^{-1/2} t_a^{1/2}
q_a^{-1/2})\}.
}
\eqnu
\label\half\eqnum*
$$

(ii) If  $\de_{2m+1}$ is given by the same formula
for $t_\nu=q_\nu^{m_\nu+1}$, then the polynomials
$\de_{2m+1} p^{(m+1)}_b$ are  anti-symmetric
 ($s_i(\  )= -(\ )$ where $1\le i\le n$)
and satisfy relations (\ref\Lf)
for $t_\nu'=q_\nu^{-m_\nu}$ and the eigenvalues equal to $f(t^{\rho}q^{b})$.
\endproclaim
The proof will be based on the following
formulas (for arbitrary $q,t$)
which are checked exactly in the same way as
the properties of the main anti-involution from
[C2] (Theorem 4.1):
$$
\eqalignno{
&(T^{q,t}_j)^\iota \ =\  \mu_{q,t} T^{q,t}_j  \mu_{q,t}^{-1},\
0\le j\le n.
 &\eqnu
\label\iotat\eqnum*
}
$$
We see that the multiplication by $\mu_{q,t}$
takes the eigenfunctions of the operators $Y_1^{q,t},\ldots,Y_n^{q,t}$
to the eigenfunctions of $\{(Y_i^{q,t})^\iota\}$ corresponding
to the same eigenvalues.
For instance, $g^{q,t}_b$ from Theorem \ref\LF calculated
for $v=\nu_R$ will go to an anti-symmetric
function $g'_b$. Indeed, the latter satisfies relations
(\ref\Lgv) for $\vep=-t^{-1}$ and $T^\iota$ instead of $T$.
It means that it is anti-symmetric.
Combining this with statement (ii) of the same
theorem we see that the multiplication
by $\mu_{q,t} d^{q,t}$ transfers the polynomials $p_b^{q, qt}$
to anti-symmetric eigenfunctions of the operators $L_f^{q,t^{-1}}$.
Here we used that the restriction of $(-s_i)$ to anti-symmetric
functions coinsides with $[\ ]_\dagger$ which is the restriction
of $s_i$ to symmetric ones ($i>0$). We can replace $t^{1/2}$ by $-t^{1/2}$
because $\{p_b\}$ depend only on $t$.

 Substituting  $tq^{-1}$
for $t$
we obtain exactly
assertion (ii). Here
 $\mu_{q,t}$ is a polynomial.

Comming to (i), we will use that for $t=q^{1/2+m}$ the
function
$$
\eqalignno{
&\mu_{q,tq^{-1}}^{-1}\prod_{a \in R_+^\vee}\prod_{i=2}^{2m}
 \{
(x_a^{1/2} t_a^{1/2}q_a^{-i/2}- x_a^{-1/2} t_a^{-1/2}q_a^{i/2})
\}
&\eqnu
\label\anti\eqnum*
}
$$
is anti-symmetric for all $s_0, s_1,\ldots, s_n$.
Hence the corresponding multiplication
turns $s_i$ into $-s_i$. Composing this with the multiplication
by $\mu_{q,tq^{-1}} d^{q,tq^{-1}}$ we get (i).

We mention that the statements can be easily extended
to the case of an arbitrary set $v$ (see Theorem \ref\LF).
Replacing $a\in R^\vee, B, r,q_a$ by $\al\in R, P, \rho, q$ one
arrives at the dual counterpart of this construction.
We also note that the above reasoning and Proposition 4.2 from [C2]
give that for arbitrary $q,t$ and $ f\in \C[x]^W$
$$
\eqalignno{
&\psi_{q,t}L_f^{q,t}\psi_{q,t}^{-1}\ =\ L_{\bar {f}} ^{q, qt^{-1}},
\ \bar {f}(x)=f(x^{-1}),&\eqnu
\cr
\psi_{q,t}\ =\
 \prod_{a \in R_+^\vee} &(x_a^{1/2}-x_a^{-1/2})^{-1}
\prod_{i=0}^\infty {(1-x_aq_a^{i}) (1-x_a^{-1}q_a^{i})
\over
(1-x_a t_a q_a^{i}) (1-x_a^{-1}t_a^{}q_a^{i})}.
\label\phi\eqnum*
}
$$

%
%
%
%      REFERENCES
%
%
%
%\vskip 15pt
\AuthorRefNames [BGG]
\references
%\medskip
%\ninerm
%\baselineskip=11pt %!
\vfil

[B]
\name{N. Bourbaki},
{ Groupes et alg\`ebres de Lie}, Ch. {\bf 4--6},
Hermann, Paris (1969).

[C1]
\bibline,
{ Integration of quantum many- body problems by affine
Knizhnik--Za\-mo\-lod\-chi\-kov equations},
Advances in Math. {106}:1 (1994)), 65--95.

[C2]
\bibline,
{ Double affine Hecke algebras and  Macdonald's
conjectures},
Annals of Mathematics {141} (1995), 191--216.

[C3]
\bibline,
{Macdonald's evaluation conjectures and difference
Fourier transform}, Inventiones Math. (1995).

[C4]
\bibline,
{Non-symmetric Macdonald polynomials},
I.M.R.N. (1995).

[He]
\name{G.J. Heckman},
{  An elementary approach to the hypergeometric shift operators of
Opdam}, Invent.Math. {  103} (1991), 341--350.

[HO]
\bibline, and \name{E.M. Opdam},
{ Root systems and hypergeometric functions I},
Comp. Math. {  64} (1987), 329--352.

[Ka]
\name {K. Kadell},
{The Dyson polynomials}, Preprint 1995.

[Ki]
\name {A. Kirillov, Jr.},
{Inner product on conformal blocks and Macdonald's
polynomials at roots of unity}, Preprint (1995).

[M1]
\name{I.G. Macdonald}, {  A new class of symmetric functions },
Publ.I.R.M.A., Strasbourg, Actes 20-e Seminaire Lotharingen,
(1988), 131--171 .

[M2]
\bibline, {  Orthogonal polynomials associated with root
systems},Preprint(1988).

[M3]
\bibline, { Affine Hecke algebras and orthogonal polynomials},
S\'eminaire Bourbaki{  47}:797 (1995), 01--18.

[O]
\name{E.M. Opdam},
{  Some applications of hypergeometric shift
operators}, Invent.Math.{  98} (1989), 1--18.

\endreferences

\bye